\documentclass[11pt,reqno]{article}
\usepackage{amsmath}
\usepackage{amssymb}

\newtheorem{theo}{Theorem}[section]\newtheorem{lem}{Lemma}[section]

\newtheorem{coro}{Corollary}[section]\numberwithin{equation}{section}
\def\pf{{\textit {Proof:} }}
\def\rm{{\textit {Remark:} }}
\def\rms{{\textit {Remarks:} }}

\newcommand{\mysection}[1]{\section{#1}\setcounter{equation}{0}}

\newfont{\bb}{msbm10 at 12pt}

\newcommand{\bal}{\begin{align}}      \newcommand{\eal}{\end{align}}
\newcommand{\ba}{\begin{array}}      \newcommand{\ea}{\end{array}}
\newcommand{\bc}{\begin{center}}     \newcommand{\ec}{\end{center}}
\newcommand{\be}{\begin{enumerate}}  \newcommand{\ee}{\end{enumerate}}
\newcommand{\beq}{\begin{eqnarray}}  \newcommand{\eeq}{\end{eqnarray}}
\newcommand{\beQ}{\begin{eqnarray*}} \newcommand{\eeQ}{\end{eqnarray*}}
\newcommand{\bi}{\begin{itemize}}    \newcommand{\ei}{\end{itemize}}
\newcommand{\bt}{\begin{tabular}}    \newcommand{\et}{\end{tabular}}
\newcommand{\bdm}{\begin{displaymath}} \newcommand{\edm}{\end{displaymath}}

\def\qed{\hfill{Q.E.D.}\smallskip}

\begin{document}

\title{Positive mass theorems for
asymptotically AdS spacetimes with arbitrary cosmological
constant}

\author{Naqing Xie \and Xiao Zhang}

\date{}

\maketitle \pagenumbering{arabic}
\begin{abstract}
We formulate and prove the Lorentzian version of the positive mass
theorems with arbitrary negative cosmological constant for
asymptotically AdS spacetimes. This work is the continuation of
the second author's recent work on the positive mass theorem on
asymptotically hyperbolic 3-manifolds.
\end{abstract}
\mysection{Introduction} In general relativity, our spacetime is
modelled by a 4-dimensional Lorentzian manifold
$(N^{1,3},\widetilde{g})$ with the Lorentzian metric
$\widetilde{g}$ of signature $(-1,1,1,1)$ which satisfies the
Einstein field equations
\begin{equation}\label{ein}
\widetilde{Ric}(\widetilde{g})-\frac{\widetilde{R}(\widetilde{g})}{2}\widetilde{g}+\Lambda
\widetilde{g}=T\end{equation} where $\widetilde{Ric}$,
$\widetilde{R}$ are the Ricci and scalar curvatures of
$\widetilde{g}$ respectively, $T$ is the energy-momentum tensor of
matter, and $\Lambda$ is the cosmological constant.

It is well known that the positive mass theorem plays an important
role in general relativity. The definition of the total energy and
the total linear momentum for asymptotically flat spaces was given
by Arnowitt-Deser-Misner from the Hamiltonian point of view
\cite{ADM}. This ADM mass is in fact a geometric invariance
\cite{B,C86}. Physicists believe, with some justification, that
the total mass for a nontrivial isolated gravitational system must
be positive. This was the famous positive mass conjecture which
was first proved by Schoen and Yau in a series of papers
\cite{SY1,SY2,SY3} using minimal surface techniques and then by
Witten \cite{Wi,PT,B} using spinors.

It is natural to extend the positive mass theorem to
asymptotically AdS spacetime where spatial infinities are
asymptotically hyperbolic. Such a theorem was proved with a fixed
negative cosmological constant for spacelike, asymptotically
hyperbolic hypersurfaces with zero second fundamental form in
\cite{Wax,CH,CN}, and with nonzero second fundamental form in
\cite{Z2,Ma}. In general, there are two versions of the positive
mass theorem (cf. \cite{Z99}). One is the Riemannian setting to
use the initial data set which is a 3-dimensional Riemannian
manifolds equipped with another 2-tensor. The other is the
Lorentzian setting to use a spacelike hypersurface in
4-dimensional Lorentzian manifolds. Although, technically,
different spin structures are used in different settings, the two
versions are essentially equivalent in asymptotically flat
spacetimes. Interestingly, the situation changes in asymptotically
AdS spacetimes and the Riemannian version of the positive mass
theorem in \cite{Z2} is not equivalent to the Lorentzian version
in \cite{Ma}. For instance, for the maximal spacelike
hypersurfaces in AdS spacetimes, the dominant energy condition in
AdS spacetimes implies the energy condition in \cite{Z2}, hence
the theorem in \cite{Z2} holds. However, this theorem is not
included in the nonnegativity of the energy-momentum matrix in
\cite{Ma}.

The present paper is essentially the continuation of the second
author's recent work in \cite{Z2}. We will prove positive mass
theorems for asymptotically AdS spacetimes with arbitrary negative
cosmological constant. We first define $e_0$- Killing spinors and
use it to obtain the corresponding Lorentzian version of the
positive mass theorem in \cite{Z2}. We then use imaginary- Killing
spinors to prove another positive mass theorem analogous to the
one in \cite{Ma}. We would like to point out that we use a little
different setting to study a spacelike hypersurface in
asymptotically AdS spacetimes in the second case, instead of
extending an initial data set to an asymptotically AdS spacetime
in \cite{Ma}. We note that it was used to study the quasi-local
mass in \cite{WY} for the positive mass theorem with zero second
fundamental form for asymptotically AdS spacetimes with arbitrary
negative cosmological constant.

It is an interesting question whether the total angular momentum
can be dominated by the total energy. In \cite{CS}, Corvino and
Schoen constructed regular solutions of vacuum Einstein constraint
equations, which are Kerr at infinity. This initial data set
indicates, in general, there is no relation between the total
energy and the total angular momentum. However, certain extra
energy conditions were found in asymptotically flat spacetimes in
\cite{Z99} that the total angular momentum is dominated by the
total energy. But the analogue of this new energy condition does
not imply the similar result in asymptotically AdS spacetimes.

This paper is organized as follows: In Section 2, we make a study
of $e_0$- as well as imaginary- Killing spinors in AdS spacetime
along the hyperbolic 3-space. Section 3 gives the definition of
total energy-momenta for asymptotically AdS initial data sets. In
Section 4, we derive a Weitzenb\"{o}ck formula for $e_0$-Killing
hypersurface Dirac-Witten operator and state some known results on
comparing two spin connections. Section 5 deals with the boundary
value problem of the Dirac-Witten equation and a positive mass
theorem is proved. In Section 6, by using the imaginary Killing
spinors, we reach another positive mass theorem which corresponds
to one of the energy-momentum inequalities from the definite
positivity of Maerten's operator.

\mysection{The AdS Spacetime, the Hyperbolic Space, and the $e_0$-
and imaginary Killing Spinors} The anti-de Sitter (AdS) spacetime
$(N^{1,3},\widetilde{g}_{AdS})$ is a static spherically solution
to the vacuum (i.e., $T=0$) Einstein equation (\ref{ein}) with
negative cosmological constant $\Lambda=\frac{-3}{l^2}$ which
reads
\begin{equation}\label{AdS1}
\widetilde{g}_{AdS}=-(\frac{\widehat{r}^2}{l^2}+1)dt^2+(\frac{\widehat{r}^2}{l^2}+1)^{-1}d\widehat{r}^2+\widehat{r}^2(d\theta^2+\sin^2\theta
d \psi^2).\end{equation} Set $\kappa=l^{-1}>0$ and
$\widehat{r}=\frac{\sinh(\kappa r)}{\kappa}$, then in terms of the
polar coordinate system $(r,\theta,\psi)$ $(0 < r <\infty$, $0
\leq \theta < \pi$ and $0 \leq \psi < 2\pi$), the AdS metric can
be rewritten as
\begin{equation}\label{AdS} \widetilde{g}_{AdS}=-\cosh^2(\kappa
r)dt^2+\overset{\circ}g\end{equation}
 where
\begin{equation}\label{hyperbolic}\overset{\circ}g=dr^2+\frac{\sinh^2(\kappa r)}{\kappa^2}(d\theta^2+\sin^2\theta
d \psi^2). \end{equation}The hyperbolic 3-space $\mathbb{H}^3$ is
the $t-$slice in $(N^{1,3},\widetilde{g}_{AdS})$ which is
topologically $\mathbb{R}^3$ endowed with the metric
$\overset{\circ}g$. Note that it is totally geodesic and has
constant sectional curvature $-\kappa^2$.

We use the convention that the Greek indices
$\alpha,\beta,\gamma,\ldots$ run over the spacetime and the Latin
ones $i,j,k,\ldots$ are the spatial indices. Denote the associated
orthonormal frame $\{\overset{\circ}e_\alpha\}$ by
$$\overset{\circ}e_0=\frac{1}{\cosh(\kappa r)}\frac{\partial}{\partial t},\ \overset{\circ}e_1=\frac{\partial}{\partial r}, \ \overset{\circ}e_2=\frac{\kappa}{\sinh(\kappa r)}\frac{\partial}{\partial \theta}, \ \overset{\circ}e_3=\frac{\kappa}{\sinh(\kappa r)\sin\theta}\frac{\partial}{\partial \psi}$$
and its coframe $\{\overset{\circ}e^\alpha\}$ by
$$\overset{\circ}e^0=\cosh(\kappa r)dt,\ \overset{\circ}e^1=dr, \ \overset{\circ}e^2=\frac{\sinh(\kappa r)}{\kappa}d\theta, \ \overset{\circ}e^3=\frac{\sinh(\kappa r)\sin\theta}{\kappa} d \psi$$
respectively.

Let $\mathbb{S}$ be the (locally) spinor bundle of
$(N^{1,3},\widetilde{g}_{AdS})$ and its restriction to
$\mathbb{H}^3$ is globally defined since every orientable
3-manifold is spin. We say that a spinor $\Phi_0 \in
\Gamma(\mathbb{S})$ is an $e_0$-Killing spinor (along
$\mathbb{H}^3$) if
\begin{equation}\label{e0k}
\nabla^{AdS}_X\Phi_0+\frac{\kappa}{2}e_0\cdot X\cdot
\Phi_0=0\end{equation} for every tangent vector $X$ of
$\mathbb{H}^3$.

Choose a standard symplectic basis as in \cite{PT} and \cite{Z01},
the spinors over AdS can be written as a 4-vector valued functions
$\Phi=(\Phi^1,\Phi^2,\Phi^3,\Phi^4)^t\in \mathbb{C}^4$. We fix the
following Clifford representation throughout the paper:
\begin{equation*}\overset{\circ}e_0 \mapsto
\begin{pmatrix}\ &\ &1 &\ \\ \ &\ & \ &1\\1&\ &\ &\ \\ \ &1&\ &\ \end{pmatrix}, \ \overset{\circ}e_1 \mapsto
\begin{pmatrix}\ &\ &-1 &\ \\ \ &\ & \ &1\\1&\ &\ &\ \\ \ &-1&\ &\ \end{pmatrix},\end{equation*}\begin{equation}\label{repre} \overset{\circ}e_2 \mapsto
\begin{pmatrix}\ &\ &\ &1 \\ \ &\ & 1 &\ \\\ &-1 &\ &\ \\ -1 &\ &\ &\ \end{pmatrix}, \ \overset{\circ}e_3 \mapsto
\begin{pmatrix}\ &\ &\ &\sqrt{-1} \\ \ &\ & -\sqrt{-1} &\ \\\ &-\sqrt{-1} &\ &\ \\ \sqrt{-1} &\ &\ &\ \end{pmatrix}.\end{equation}
Using it, one has \begin{lem} The set of the solutions of the
$e_0$-Killing equation (\ref{e0k}) is 4-dimensional. Precisely,
\begin{equation}\label{KS}\Phi_0=\begin{pmatrix}(\lambda_1e^{\frac{\sqrt{-1}}{2}\psi}\sin\frac{\theta}{2}+\lambda_2e^{-\frac{\sqrt{-1}}{2}\psi}\cos\frac{\theta}{2})e^{-\frac{\kappa r}{2}}\\(\lambda_2e^{-\frac{\sqrt{-1}}{2}\psi}\sin\frac{\theta}{2}-\lambda_1e^{\frac{\sqrt{-1}}{2}\psi}\cos\frac{\theta}{2})e^{\frac{\kappa r}{2}}\\(\lambda_3e^{\frac{\sqrt{-1}}{2}\psi}\sin\frac{\theta}{2}+\lambda_4e^{-\frac{\sqrt{-1}}{2}\psi}\cos\frac{\theta}{2})e^{\frac{\kappa r}{2}}\\(\lambda_4e^{-\frac{\sqrt{-1}}{2}\psi}\sin\frac{\theta}{2}-\lambda_3e^{\frac{\sqrt{-1}}{2}\psi}\cos\frac{\theta}{2})e^{-\frac{\kappa r}{2}}\end{pmatrix}\in \mathbb{C}^4\end{equation}
where $\lambda_1$, $\lambda_2$, $\lambda_3$, and $\lambda_4$ are
four arbitrary complex numbers.
\end{lem}
Due to the fact that $e_0$ changes the chirality of spinors, the
form of $e_0$-Killing spinors looks different from that of the
imaginary Killing spinors \cite{L,HT}. In fact,
\begin{lem}The
imaginary Killing spinors along $\mathbb{H}^3$ satisfying the
imaginary Killing equations
\begin{equation}\label{ka}
\nabla_X \Phi+\frac{\kappa \sqrt{-1}}{2}X\cdot\Phi=0\ \ \mbox{for
each tangent vector $X$ of $\mathbb{H}^3$}\end{equation} are of
the form \begin{equation}\label{imk}\begin{pmatrix}\Phi^1\\\Phi^2\\
\Phi^3\\\Phi^4\end{pmatrix}=\begin{pmatrix}u^+e^{\frac{\kappa
r}{2}}+u^-e^{-\frac{\kappa r}{2}}\\v^+e^{\frac{\kappa
r}{2}}+v^-e^{-\frac{\kappa r}{2}}\\-\sqrt{-1}u^+e^{\frac{\kappa
r}{2}}+\sqrt{-1}u^-e^{-\frac{\kappa r}{2}}
\\
\sqrt{-1}v^+e^{\frac{\kappa r}{2}}-\sqrt{-1}v^-e^{-\frac{\kappa
r}{2}}\end{pmatrix}\end{equation} where
\begin{eqnarray*}
u^+ &=& \lambda_1
e^{\frac{\sqrt{-1}}{2}\psi}\sin\frac{\theta}{2}+\lambda_2
e^{\frac{-\sqrt{-1}}{2}\psi}\cos\frac{\theta}{2},\\
u^-&=&\lambda_3
e^{\frac{\sqrt{-1}}{2}\psi}\sin\frac{\theta}{2}+\lambda_4
e^{\frac{-\sqrt{-1}}{2}\psi}\cos\frac{\theta}{2},\\
v^+&=&-\lambda_3
e^{\frac{\sqrt{-1}}{2}\psi}\sin\frac{\theta}{2}+\lambda_4
e^{\frac{-\sqrt{-1}}{2}\psi}\cos\frac{\theta}{2}, \\
v^-&=&-\lambda_1
e^{\frac{\sqrt{-1}}{2}\psi}\sin\frac{\theta}{2}+\lambda_2
e^{\frac{-\sqrt{-1}}{2}\psi}\cos\frac{\theta}{2}.\end{eqnarray*}
Here $\lambda_1$, $\lambda_2$, $\lambda_3$, and $\lambda_4$ are
four arbitrary complex numbers.\end{lem}
\rm \\
The imaginary Killing spinors in the full spacetime look like in
the similar form but $\lambda_\mu$ is time dependent \cite{HT},
i.e.
$$\lambda_1=C_1\cos (\frac{\kappa}{2}t)+C_3\sin(\frac{\kappa}{2}
t), \ \ \lambda_3=C_3\cos(\frac{\kappa}{2} t)- C_1\sin
(\frac{\kappa}{2}t)$$ $$\lambda_2=C_2\cos(\frac{\kappa}{2}
t)+C_4\sin (\frac{\kappa}{2}t), \ \ \lambda_4=C_4\cos
(\frac{\kappa}{2}t)- C_2\sin(\frac{\kappa}{2} t)$$ where $C_1,\
C_2,\ C_3,\ C_4$ are four arbitrary complex constants.

\mysection{Definition of the Total Energy-Momenta} In this
section, we define the total energy-momenta for asymptotically AdS
initial data sets. Suppose that $(N^{1,3},\widetilde{g})$ is a
Lorentzian manifold with the Lorentzian metric $\widetilde{g}$ of
signature $(-1,1,1,1)$ satisfying the Einstein field equations.
Usually, a triple $(M,g_{ij},h_{ij})$ is served as a Cauchy
surface on the initial problem of the Einstein equations. Here $M$
is a 3-dimensional spacelike hypersurface with induced Riemannian
metric $g_{ij}$ and $h_{ij}$ is the second fundamental form of $M$
in $N$. We say that the initial data set $(M,g_{ij},h_{ij})$ is
asymptotically AdS if:\\
(1) There is a compact set $K  \subset M$ such that $M_\infty=M
\setminus K$ is diffeomorphic to $\mathbb{R}^3 \setminus
\mbox{open ball}$; \\
(2) Under this diffeomorphism, the metric
$g_{ij}=g(\overset{\circ}e_i,\overset{\circ}e_j)$ on the end
$M_\infty$ is of the form
$$g_{ij}=\delta_{ij}+a_{ij}$$
where $a_{ij}$ satisfies
\begin{equation}\label{asymg}
a_{ij}=O(e^{- \tau \kappa r}), \
\overset{\circ}{\overline{\nabla}}_k a_{ij}=O(e^{- \tau \kappa
r}), \
\overset{\circ}{\overline{\nabla}}_l\overset{\circ}{\overline{\nabla}}_k
a_{ij}=O(e^{-\tau \kappa r});\end{equation} and the second
fundamental form $h_{ij}=h(\overset{\circ}e_i,\overset{\circ}e_j)$
satisfies
\begin{equation}\label{asymh}
h_{ij}=O(e^{- \tau \kappa r}), \
\overset{\circ}{\overline{\nabla}}_k h_{ij}=O(e^{-\tau \kappa
r})\end{equation} for $\tau >\frac{3}{2}$. Here
$\overset{\circ}{\overline{\nabla}}$ is the Levi-Civita connection
with respect to
the hyperbolic metric $\overset{\circ}g$;\\
(3) There exists a distance function $\rho _z$ such that
\begin{equation}T _{00} e ^{\rho _z} \in L ^1 (M),\;\;\;\;\;\;T _{0i} e
^{\rho _z }\in L ^1 (M).\end{equation}

Denote $n^0=1$, $n^i\ (i=1,2,3)$ the restriction of the natural
coordinate $x^i$ to the unit round sphere, i.e.
$$n^0=1, \ n^1=\sin\theta\cos\psi, \ n^2=\sin\theta\sin\psi, \ n^3=\cos \theta$$
and
$$\mathcal{\varepsilon}_i=\overset{\circ}{{\overline{\nabla}}^j
}g_{ij}-\overset{\circ}
{\overline{\nabla}}_itr_{\overset{\circ}g}(g)-\kappa(a_{1i}-g_{1i}tr_{\overset{\circ}g}(a)),$$
$$\mathcal{P}_{ki}=h_{ki}-g_{ki}tr_{\overset{\circ}g}(h).$$

For such a spacetime, (under a fixed diffeomorphism) the total
energy vector $E_{\{\nu\}}$ and for each $k$ the total linear
momentum vector $P_{\{\nu\}k}$ are defined by
\begin{equation}
E_{\{\nu\}}=\frac{1}{16\pi}\lim_{r\rightarrow
\infty}\int_{S_r}\mathcal{\varepsilon}_1\omega_\nu,
\end{equation}
\begin{equation}
P_{\{\nu\}k}=\frac{1}{8\pi}\lim_{r\rightarrow
\infty}\int_{S_r}\mathcal{P}_{k1}\omega_\nu
\end{equation}
where\begin{equation}\label{omega} \omega_\nu=n^\nu e^{\kappa
r}\overset{\circ}e^2\wedge \overset{\circ}e^3 \mbox{
$(\nu=0,1,2,3)$.}\end{equation}\rms \\
(1) For simplicity, we just assume that there is only one end. The
extension of
multi-ends case is straightforward.\\
(2) Our definition $E_{\{\nu\}}$ is the same as $p_{\nu}$ in
\cite{CH}. The geometric invariance of the total energy was given
in \cite{CN,CH} as well as in \cite{Wax} for the case with a
spherical conformal infinity.\\
(3) The definition of the total linear momentum vector can be
found in \cite{Z2}. In fact, the Lorentzian lengths of
$P_{\{\nu\}1}$ is invariant but $P_{\{\nu\}A}$ ($A=2,3$) is not
\cite[Prof. 2.1]{Z2}. Therefore, for $i=1,2,3$,
\begin{equation}
(c_1E_{\{0\}}+c_2P_{\{0\}1})^2-\sum_{i}(c_1E_{\{i\}}+c_2P_{\{i\}1})^2\end{equation}
gives a geometric invariant where $c_1$ and $c_2$ are real
constants.

\mysection{The Spin Connections, the Dirac-Witten Operators, and
the Weitzenb\"{o}ck Formula} In this section, we establish a
Weitzenb\"{o}ck formula for the Dirac-Witten operator associated
with the $e_0$-Killing connection. We also state some known
results (due to Min-Oo \cite{M} and certain generalizations in
\cite{AnD,He3,Z1}) on comparing two spin connections.

Recall that $(N^{1,3},\widetilde{g})$ is a Lorentzian manifold
with the Lorentzian metric $\widetilde{g}$ of signature
$(-1,1,1,1)$ satisfying the Einstein field equations. Let
$(M,g,h)$ be a 3-dimensional spacelike hypersurface with induced
Riemannian metric $g_{ij}$ and $h_{ij}$ is the second fundamental
form of $M$ in $N$. Let $\mathbb{S}$ be the (locally) spinor
bundle of $N$ and we still denote by $\mathbb{S}$ its restriction
to $M$. Let $\nabla$ and $\overline{\nabla}$ be the Levi-Civita
connections of $\widetilde{g}$ and $g$ respectively. We also
denote by the same symbols their lifts to the spinor bundle
$\mathbb{S}$.

Fix a point $p\in M$ and an orthonormal basis $\{e_\alpha\}$ of
$T_pN$ with $e_0$ normal and $\{e_i\}$ tangent to $M$. Extend
$\{e_\alpha\}$ to a local orthonormal frame in a neighborhood of
$p$ in $M$ such that $(\overline{\nabla}_ie_j)_p=0$. Extend this
to a local orthonormal frame $\{e_\alpha\}$ for $N$ with
$(\nabla_0e_j)_p=0$. Let $\{e^\alpha\}$ be its dual frame. Then
\begin{equation}(\nabla_ie_j)_p=h_{ij}e_0, \
(\nabla_ie_0)_p=h_{ij}e_j\end{equation} where
$h_{ij}=\widetilde{g}(\nabla_ie_0,e_j)$ are the components of its
second fundamental form at $p$. The two connections on the spinor
bundle are related by
\begin{equation}
\nabla_i=\overline{\nabla}_i-\frac{1}{2}h_{ij}e_0\cdot
e_j\cdot.\end{equation} We define the $e_0$-Killing connection by
\begin{equation}
\overset{\wedge}\nabla_X=\nabla_X+\frac{\kappa}{2}e_0\cdot X\cdot.
\end{equation}
Then the associated hypersurface Dirac-Witten operators are
\begin{eqnarray}
D=\sum_{k=1}^3e_k\cdot\nabla_k,
\end{eqnarray}
\begin{eqnarray}
\overset{\wedge}D=\sum_{k=1}^3e_k\cdot\overset{\wedge}\nabla_k=D+\frac{3\kappa}{2}e_0\cdot\end{eqnarray}
There are two choices of metrics on the spinor bundle
\cite{PT,Z01}. Restricted $\mathbb{S}$ to M inherits an Hermitian
metric $(\phi,\psi)$ and a positive definite metric $<\phi,\psi>$.
They are related by the equation
$$(\phi,\psi)=<e_0\cdot\phi,\psi>.$$
With respect to $<\cdot,\cdot>$, $e_i$ is skew-Hermitian while
$e_0$ is Hermitian \cite{PT,Z01}. We also note that
$\overline{\nabla}$ is compatible with $<\cdot,\cdot>$ but
$\nabla$ is not. Moreover,
\begin{equation}
\overline{\nabla}_i(e_0\cdot\phi)=e_0\cdot \overline{\nabla}_i
\phi.
\end{equation}

As usual, we have the following formulae.
\begin{eqnarray}
& &d(<\phi, \nabla_i\psi>\ast e^i) =(e_i<\phi,\nabla_i\psi>)\ast 1  \nonumber\\
& = &((\nabla_ie_0)\cdot\phi, \nabla_i\psi)\ast 1+<\nabla_i\phi,
\nabla_i\psi>\ast 1+ <\phi, \nabla_i\nabla_i\psi>\ast
1    \nonumber\\
&=&(h_{ij}e_j\cdot \phi, \nabla_i\psi)\ast 1+<\nabla_i\phi,
\nabla_i\psi>\ast 1+ <\phi, \nabla_i\nabla_i\psi>\ast1  \nonumber\\
&=&<h_{ij}e_0\cdot e_j\cdot \phi, \nabla_i\psi>\ast
1+<\nabla_i\phi, \nabla_i\psi>\ast 1+ <\phi,
\nabla_i\nabla_i\psi>\ast1
\end{eqnarray}
and
 \begin{eqnarray}
  d(<e_i\cdot\phi,\psi>\ast
e^i)&=&(<D\phi,\psi>-<\phi,D\psi>)\ast1   \\
&=&(<\overset{\wedge}D\phi,\psi>-<\phi,\overset{\wedge}D\psi>)\ast1.
 \end{eqnarray}
Hence,
\begin{equation}
\nabla_i^\ast=-\nabla_i-h_{ij}e_0\cdot e_j\cdot
\end{equation}
\begin{equation}
\overset{\wedge}\nabla^\ast_i=\nabla_i^\ast+\frac{\kappa}{2}e_0\cdot
e_i
\end{equation}
\begin{equation}
D^\ast=D
\end{equation}
\begin{equation}
\overset{\wedge}D^\ast=\overset{\wedge}D=D+\frac{3\kappa}{2}e_0\cdot=D^\ast+\frac{3\kappa}{2}e_0\cdot\end{equation}with
respect to $<\cdot, \cdot>$.

The corresponding Weitzenb\"{o}ck formula for the Dirac-Witten
operator $\overset{\wedge}D$ is therefore:
\begin{lem}One has
\begin{equation}\label{W}\overset{\wedge}D^\ast\overset{\wedge}D=\overset{\wedge}\nabla^\ast\overset{\wedge}\nabla+
\overset{\wedge}{\mathcal{R}}
\end{equation} where
$$\overset{\wedge}{\mathcal{R}}=\frac{1}{4}(Scal^{\widetilde{g}}+2\widetilde{R}_{00}+2\widetilde{R}_{0i}e_0\cdot e_i\cdot+6\kappa^2-4\kappa tr(h)) \in
End(\mathbb{S}).$$\end{lem} \pf By straightforward computation, it
follows that
\begin{eqnarray*}
\overset{\wedge}D^\ast\overset{\wedge}D&=&\overset{\wedge}D^2=(D+\frac{3\kappa}{2}e_0\cdot)\circ(D+\frac{3\kappa}{2}e_0\cdot)\\
&=&D^2+\frac{9\kappa^2}{4}+\frac{3\kappa}{2}e_k\cdot(\nabla_ke_0)\\
&=&D^2+\frac{9\kappa^2}{4}+\frac{3\kappa}{2}e_k\cdot(h_{kj}e_j\cdot)\\
&=&D^2+\frac{9\kappa^2}{4}-\frac{3\kappa}{2}tr(h).\end{eqnarray*}
On the other hand, \begin{eqnarray*}\overset{\wedge}{\nabla}^\ast
\overset{\wedge}{\nabla}&=&(\nabla_i^\ast+\frac{\kappa}{2}e_0\cdot
e_i\cdot)\circ(\nabla_i+\frac{\kappa}{2}e_0\cdot e_i\cdot)\\
&=&(-\nabla_i-h_{ij}e_0\cdot e_j\cdot+\frac{\kappa}{2}e_0\cdot
e_i\cdot)\circ(\nabla_i+\frac{\kappa}{2}e_0\cdot e_i\cdot)\\
&=&\nabla^\ast\nabla-\frac{\kappa}{2}(\nabla_ie_0)\cdot
e_i\cdot-\frac{\kappa}{2}e_0\cdot
(\nabla_ie_i)\cdot-\frac{\kappa}{2}h_{ij}e_0\cdot e_j\cdot
e_0\cdot e_i+\frac{3\kappa^2}{4}\\
&=&\nabla^\ast\nabla+\frac{3\kappa^2}{4}-\frac{\kappa}{2}tr(h).\end{eqnarray*}
The standard Weitzenb\"{o}ck formula \cite{Wi,PT} reads
$$D^2=\nabla^\ast\nabla+\frac{1}{4}(Scal^{\widetilde{g}}+2\widetilde{R}_{00}+2\widetilde{R}_{0i}e_0\cdot
e_i\cdot).$$This proves the lemma.\qed

Since we work on a non-compact manifold we need the following
integrated version of the Weitzenb\"{o}ck formula (\ref{W})
involving a boundary term:
\begin{lem}
We have \begin{equation}\label{WI}\int_{\partial M}<\phi,
\overset{\wedge}\nabla_i \phi+e_i\cdot\overset{\wedge}D\phi>\ast
e^i=\int_M\{|\overset{\wedge}\nabla\phi|^2-|\overset{\wedge}D\phi|^2\}\ast
1
+\int_M<\phi,\overset{\wedge}{\mathcal{R}}\phi>\ast1\end{equation}for
all $\phi \in \Gamma(\mathbb{S})$.
\end{lem}
The assumption we make in order to prove the positive mass theorem
is a modified version of dominant energy condition:
\begin{equation}\label{energy}
Scal^{\widetilde{g}}+2\widetilde{R}_{00}+6\kappa^2-4\kappa
tr(h)¡¡\geq \sqrt{\sum_i (2\widetilde{R}_{0i})^2}.\end{equation}
This ensures
\begin{equation}\label{2energy}
<\phi,\overset{\wedge}{\mathcal{R}}\phi>\geq 0 ,\ \forall \phi \in
\Gamma(\mathbb{S}).
\end{equation}
\rms \\
(1) Set $\kappa=1$, $p_{ij}=-h_{ij}+\delta_{ij}$. Let $$\mu
=\frac{1}{2}(Scal^g+(p^i_i)^2-p_{ij}p^{ij}),$$
$$\overline{\omega}_j
=\overline{\nabla}^ip_{ji}-\overline{\nabla}_jp^i_i.$$ Then by the
Gauss-Codazzi equations of $M$ in $N$,  one can show that the
energy condition (\ref{energy}) is equivalent to the following
``dominant energy condition" in \cite{Z2}
$$\mu \geq
\sqrt{\sum_i\overline{\omega}_i^2}.$$ \\ (2) If $M$ is maximal,
i.e. $tr(h)=0$, then the energy condition (\ref{energy}) just
reduces to the standard dominant energy condition
\begin{equation}\label{sdec}
T_{00}\geq \sqrt{\sum_iT_{0i}^2}.\end{equation} Einstein equation
(\ref{ein}) gives
$$Scal^{\widetilde{g}}=2(\widetilde{G}_{00}-\widetilde{R}_{00})$$
and
$$\widetilde{R}_{0i}e_0\cdot e_i \cdot =e_0\cdot (\widetilde{G}_{0i}e_i)\cdot$$
where
$\widetilde{G}_{\mu\nu}=\widetilde{R}_{\mu\nu}-\frac{\widetilde{R}}{2}\widetilde{g}_{\mu\nu}$
is the Einstein tensor. Then
\begin{eqnarray*}\overset{\wedge}{\mathcal{R}}&=&(\frac{1}{2}\widetilde{G}_{00}e_0\cdot-\frac{3\kappa^2}{2}\widetilde{g}_{00}e_0\cdot-\frac{1}{2}\widetilde{G}_{0i}e_i\cdot)e_0\cdot\\
&=&\frac{1}{2}(T_{00}e_0\cdot-T_{0i}e_i\cdot)e_0\cdot.\end{eqnarray*}\\[0.3cm]

In order to compute the boundary term which gives rise to the
total energy and the total momentum, we define a new connection
and see the difference of two connections on the spinor bundle.
Most of the results here are due to the work in \cite{M,AnD,Z1}.
Recall that $g=\overset{\circ}g+a$ with $a=O(e^{-\tau\kappa r}), \
\overset{\circ}{\overline{\nabla}} a=O(e^{-\tau\kappa r}), \
\overset{\circ}{\overline{\nabla}}\overset{\circ}{\overline{\nabla}}a=O(e^{-\tau\kappa
r})$. Orthonormalizing $\overset{\circ}e_i$ with respect to
$\overset{\circ}g$ gives rise an orthonormal basis $e_i$ with
respect to $g$, i.e.
\begin{equation}
e_i=\overset{\circ}e_i-\frac{1}{2}a_{ik}\overset{\circ}e_k+o(e^{-\tau\kappa
r}).\end{equation} This gives a gauge transformation
$$\mathcal{A}: SO(\overset{\circ}g) \rightarrow SO(g)$$
$$\overset{\circ}e_i \mapsto e_i$$ (and in addition $e_0 \mapsto e_0$)
which identifies the corresponding spin group and spinor bundles.

To compare $\overline{\nabla}$ and
$\overset{\circ}{\overline{\nabla}}$ in particular their lifts to
the spinor bundles, one introduces a new connection
$\widetilde{\nabla}=\mathcal{A}\circ\overset{\circ}{\overline{\nabla}}\circ\mathcal{A}^{-1}$.
This new connection is compatible with the metric $g$ but has a
torsion
\begin{eqnarray}
\widetilde{T}(X,Y)&=&\widetilde{\nabla}_XY-\widetilde{\nabla}_YX-[X,Y] \nonumber\\
&=&-(\overset{\circ}{\overline{\nabla}}_X\mathcal{A})\mathcal{A}^{-1}Y+(\overset{\circ}{\overline{\nabla}}_Y\mathcal{A})\mathcal{A}^{-1}X.\label{T1}
\end{eqnarray}
Then the difference of $\overline{\nabla}$ of $\widetilde{\nabla}$
is then expressible in terms of the torsion
\begin{equation}\label{T2}
2g(\widetilde{\nabla}_XY-\overline{\nabla}_XY,Z)=g(\widetilde{T}(X,Y),Z)-g(\widetilde{T}(X,Z),Y)-g(\widetilde{T}(Y,Z),X)\end{equation}
for any tangent vectors $X,Y,Z \in TM$.

Since both $\overline{\nabla}$ and $\widetilde{\nabla}$ are
$g$-compatible, their induced connections on the spinor bundle
$\mathbb{S}(M)$ differ by
\begin{equation}
\overline{\nabla}_j-\widetilde{\nabla}_j=-\frac{1}{4}\sum_{k,l}(\omega_{kl}(e_j)-\widetilde{\omega}_{kl}(e_j))e_k\cdot
e_l\cdot
\end{equation}
where $\omega_{kl}(e_j)=-g(\overline{\nabla}_je_k,e_l)$ and
$\widetilde{\omega}_{kl}(e_j)=-g(\widetilde{\nabla}_je_k,e_l)$.

From (\ref{T1}) and (\ref{T2}) we have obtained the following
asymptotic formula
\begin{equation}
\overline{\nabla}_j-\widetilde{\nabla}_j=\frac{1}{8}\sum_{k\neq
l}(\overset{\circ}{{\overline{\nabla}}^k
}g_{jl}-\overset{\circ}{{\overline{\nabla}}^l }g_{jk})e_k\cdot
e_l\cdot+o(e^{-\tau\kappa r})\end{equation} for the difference of
the two connections acting on spinors. And further we have
\begin{lem} (Prop. 3.2, \cite{Z2}) Let $(M,g_{ij},h_{ij})$ be a 3-dimensional asymptotically AdS initial data set. Then
\begin{equation}
\sum_{j,\ j\neq i}Re<\phi,e_i\cdot
e_j\cdot(\overline{\nabla}_j-\widetilde{\nabla}_j)\phi>=\frac{1}{4}(\overset{\circ}{{\overline{\nabla}}^j}
g_{ij}-\overset{\circ}{\overline{\nabla}}_i
tr_{\overset{\circ}g}(g)+o(e^{-\tau\kappa r}))|\phi|^2
\end{equation}
for all $\phi \in \Gamma (\mathbb{S})$.\end{lem}

We extend the $e_0$-Killing spinors $\Phi_0$ in (\ref{KS}) on the
end to the inside smoothly. With respect to the metric $g$, these
$e_0$-Killing spinors $\Phi_0$ can be written as
$\overline{\Phi}_0=\mathcal{A}\Phi_0$. Let
$\overset{\wedge}{\widetilde{\nabla}}_X=\widetilde{\nabla}_X+\frac{\kappa}{2}e_0\cdot
X\cdot $. Then
\begin{eqnarray*}
\overset{\wedge}{\widetilde{\nabla}}_j\overline{\Phi}_0&=&\mathcal{A}({\overset{\circ}{\overline{\nabla}}}_j\Phi_0)+\frac{\kappa}{2}e_0\cdot
e_j\cdot(\mathcal{A}\Phi_0)\\
&=&\frac{\kappa}{4}a_{jk}e_0\cdot(\mathcal{A}\overset{\circ}e_k)\cdot\overline{\Phi}_0+o(e^{-\tau\kappa
r})\overline{\Phi}_0.\end{eqnarray*}

\mysection{The Dirac-Witten Equation and Positive Mass Theorem I} As
explained in the introduction, in this section, we will study an
elliptic boundary value problem on $M$ with given boundary values as
$r \rightarrow \infty$. We will solve for a spinor $\phi$ satisfying
the first order elliptic Dirac-Witten equation
$\overset{\wedge}D\phi=0$ on the manifold $M$ which is asymptotic to
the $e_0$-Killing spinor $\Phi_0$ at infinity. Our positive mass
theorem is then a consequence of the nice Weitzenb\"{o}ck formula
(\ref{WI}).

Let $C_0^\infty(\mathbb{S})$ be the space of smooth sections of
the spinor bundle $ \mathbb{S}$ with compact support. Define an
inner product on $\mathbb{S}$ by
\begin{equation}
(\phi,\psi)_1=\int_M
\Big{\{}<\nabla\phi,\nabla\psi>+\frac{3\kappa^2}{4}<\phi,\psi>\Big{\}}\ast
1.
\end{equation}
Let $H^1(\mathbb{S})$ be the closure of $C_0^\infty(\mathbb{S})$
with respect to this inner product. Then $H^1(\mathbb{S})$ with
the above inner product is a Hilbert space. Now define a bounded
bilinear form $\mathcal{B}$ on $C_0^\infty(\mathbb{S})$ by
\begin{equation}
\mathcal{B}(\phi,\psi)=\int_M<\overset{\wedge}D\phi,
\overset{\wedge}D\psi
>\ast 1.\end{equation}
By the Weitzenb\"{o}ck formula (\ref{WI}), we have
\begin{equation}\mathcal{B}(\phi,\phi)=\int_M
|\overset{\wedge}\nabla\phi|^2 \ast 1+\int_M<\phi,
\overset{\wedge}{\mathcal{R}}\phi>\ast 1.
\end{equation}
Due to the energy condition (\ref{energy}), we can extend
$\mathcal{B}(\cdot,\cdot)$ to $H^1(\mathbb{S})$ as a coercive (not
strictly coercive in general) bilinear form.

Take $\Phi_0$ as the $e_0$-Killing spinor in (\ref{KS}). The same
as in \cite{Z2}, due to the asymptotic conditions (\ref{asymg})
and (\ref{asymh}), we know that $\overset{\wedge}\nabla
\overline{\Phi}_0 \in L^2(\mathbb{S})$ and hence
$\overset{\wedge}D \overline{\Phi}_0 \in L^2(\mathbb{S})$. Note
that $\overline{\Phi}_0$ itself is not in $L^2(\mathbb{S})$ since
$|\overline{\Phi}_0|^2= O(e^{\kappa r})$.
\begin{lem}\label{DW}Let $(M,g_{ij},h_{ij})$ be a 3-dimensional asymptotically AdS initial data set which satisfies the energy condition (\ref{energy}). There exists a unique spinor $\Phi_1$ in $H^1(\mathbb{S})$  such that
\begin{equation}
\overset{\wedge}D(\Phi_1+\overline{\Phi}_0)=0.\end{equation}
\end{lem}
\pf Here we follow the argument of Lemma 4.2 in \cite{Z2}. Since
$\mathcal{B}(\cdot,\cdot)$ is coercive on $H^1(\mathbb{S})$, and
$\overset{\wedge}D \overline{\Phi}_0 \in L^2(\mathbb{S})$,
$\overset{\wedge}\nabla \overline{\Phi}_0 \in L^2(\mathbb{S})$,
thanks to the Lax-Milgram theorem (cf. Theorem 7.21 \cite{F}),
there exists a spinor $\Phi_1 \in H^1(\mathbb{S})$ such that
$$\overset{\wedge}D^\ast\overset{\wedge}D\Phi_1=-\overset{\wedge}D^\ast \overset{\wedge}D\overline{\Phi}_0$$
weakly. Let $\phi=\Phi_1+\overline{\Phi}_0$ and
$\psi=\overset{\wedge}D\phi$. The elliptic regularity tells us
that $\psi \in H^1(\mathbb{S})$, and
$$\overset{\wedge}D^\ast\psi=0$$
in the classical sense. Then (\ref{WI}) implies that
$\overset{\wedge}\nabla \psi=0$. We thus have $|\partial_i \log
|\psi|^2| \leq (\kappa+|h|)$ on the complement of the zero set of
$\psi$ on $M$. If there exists $x_0 \in M$ such that
$|\psi(x_0)|\neq 0$, then integrating it along a path from $x_0\in
M$ gives
$$|\psi(x)|^2 \geq |\psi(x_0)|^2e^{(\kappa+|h|)(|x_0|-|x|)}.$$
Obviously, $\psi$ is not in $L^2(\mathbb{S})$ which gives the
contradiction. Hence $\psi=0$, and the proof of this lemma is
complete.\qed

Now we state our first positive mass theorem.
\begin{theo}\label{PET} Let $(M,g_{ij},h_{ij})$ be a 3-dimensional asymptotically AdS initial data set which satisfies the energy condition (\ref{energy}).
Then the following $4 \times 4$ Hermitian matrix
\begin{equation}\label{Q1}
\begin{pmatrix}
E_{\{0\}}+P_{\{0\}1}& -E_{\{1\}}-P_{\{1\}1}& \ & \ \\
+E_{\{3\}}+P_{\{3\}1}& +\sqrt{-1}(E_{\{2\}}+P_{\{2\}1})& \ & \ \\
\ & \ & \ & \ \\
-E_{\{1\}}-P_{\{1\}1}& E_{\{0\}}+P_{\{0\}1}& \ & \ \\
-\sqrt{-1}(E_{\{2\}}+P_{\{2\}1})& -E_{\{3\}}-P_{\{3\}1}& \ & \ \\
\ & \ & \ & \ \\
\ & \ & E_{\{0\}}+P_{\{0\}1} & E_{\{1\}}+P_{\{1\}1}\\
\ & \ & -E_{\{3\}}-P_{\{3\}1} & +\sqrt{-1}(E_{\{2\}}+P_{\{2\}1})\\
\ & \ & \ & \ \\
\ & \ & E_{\{1\}}+P_{\{1\}1} & E_{\{0\}}+P_{\{0\}1}\\
\ & \ & -\sqrt{-1}(E_{\{2\}}+P_{\{2\}1})& E_{\{3\}}+P_{\{3\}1}

\end{pmatrix}\end{equation}
is positive definite. Moreover, if $E_{\{0\}}+P_{\{0\}1} =0$, then
the following equations hold on $M$:
\begin{equation}\label{gauss}R_{ijkl}+\widetilde{h}_{ik}\widetilde{h}_{jl}-\widetilde{h}_{il}\widetilde{h}_{jk}=0,\end{equation}\begin{equation}\label{codazzi}\overline{\nabla}_i\widetilde{h}_{jk}-\overline{\nabla}_j\widetilde{h}_{ik}=0\end{equation}where
$R_{ijkl}$ is the Riemann curvature tensor of $(M,g)$ and
$\widetilde{h}_{ij}=-h_{ij}+\kappa\delta_{ij}$.
\end{theo}

The positivity of the $2\times 2$ principal minor in (\ref{Q1})
also implies
\begin{coro}
In particular, we have \begin{equation}E_{\{0\}}+P_ {\{0\}1}\geq
\sqrt{\sum^3_{i=1}(E_{\{i\}}+P_
{\{i\}1})^2}.\end{equation}\end{coro}

As mentioned in the previous section, when $M$ is maximal, the
energy condition reduces to the standard dominant energy condition
and can be expressed in terms of the energy-momentum tensor
$T_{\mu\nu}$. One thus also has the result below.
\begin{coro}
Let $(M,g_{ij},h_{ij})$ be a 3-dimensional asymptotically AdS
initial data set. Assume that $M$ is maximal and satisfies the
dominant energy condition (\ref{sdec}), then the $4 \times 4$
Hermitian matrix in (\ref{Q1}) is positive definite.\end{coro}

Now we are going to prove Theorem \ref{PET}.

\pf Let $\phi$ be the solution of the Dirac-Witten equation
$\overset{\wedge}D\phi=0$ as in Lemma \ref{DW}. Submitting this
$\phi$ into the Weitzenb\"{o}ck formula (\ref{WI}), we obtain that
the boundary term is nonnegative due to the energy condition
(\ref{energy}).

Denote
 \begin{eqnarray*}
\alpha_i&=&\overset{\circ}{{\overline{\nabla}}^j}
g_{ij}-\overset{\circ}{\overline{\nabla}}_i
tr_{\overset{\circ}g}(g),\\
b_{ij}&=&-2h_{ij}+\kappa a_{ij},\\
\tau_{ki}&=&b_{ki}-g_{ki}tr_{\overset{\circ}g}(b),
 \end{eqnarray*}
and
 \begin{eqnarray*}
\beta_\nu=\frac{1}{16\pi}\lim_{r \rightarrow
\infty}\int_{S_r}(\alpha_1-\tau_{11})\omega_\nu
 \end{eqnarray*}
where $\omega_\nu$ is defined in (\ref{omega}).

Therefore we obtain
 \begin{eqnarray*}
&&\int_M|\overset{\wedge}\nabla\phi|^2\ast 1
+\int_M<\phi,\overset{\wedge}{\mathcal{R}}\phi>\ast1\\
&=&\lim_{r\rightarrow
\infty}Re\int_{S_r}<\overline{\Phi}_0,\sum_{i,j,i\neq j }e_i\cdot
e_j\cdot\overset{\wedge}\nabla_j \overline{\Phi}_0>\ast e^i\\
&=&\lim_{r\rightarrow
\infty}Re\int_{S_r}<\overline{\Phi}_0,\sum_{i,j,i\neq j }e_i\cdot
e_j\cdot(\overline{\nabla}_j-
\widetilde{\nabla}_j)\overline{\Phi}_0>\ast e^i\\
&\ &+\lim_{r\rightarrow
\infty}Re\int_{S_r}<\overline{\Phi}_0,\sum_{i,j,i\neq j }e_i\cdot
e_j\cdot\overset{\wedge}{
\widetilde{\nabla}}_j)\overline{\Phi}_0>\ast e^i\\
&\ &-\lim_{r\rightarrow
\infty}Re\int_{S_r}<\overline{\Phi}_0,\sum_{i,j,i\neq j
}\frac{1}{2}h_{jk} e_i\cdot e_j\cdot e_0\cdot e_k\cdot
\overline{\Phi}_0>\ast e^i\\
&=&\frac{1}{4}\lim_{r\rightarrow
\infty}\int_{S_r}(\overset{\circ}{{\overline{\nabla}}^j}
g_{1j}-\overset{\circ}{\overline{\nabla}}_1
tr_{\overset{\circ}g}(g))|\Phi_0|^2\overset{\circ}{\ast}\overset{\circ}{e}^1\\
&\ &+\frac{1}{4}\lim_{r\rightarrow
\infty}\int_{S_r}\kappa(a_{k1}-g_{k1}tr_{\overset{\circ}g}(a))<\Phi_0,\overset{\circ}{e}_0\cdot\overset{\circ}{e}_k\cdot\Phi_0
>\overset{\circ}{\ast}\overset{\circ}{e}^1\\
&\ &-\frac{1}{2}\lim_{r\rightarrow
\infty}\int_{S_r}(h_{k1}-g_{k1}tr_{\overset{\circ}g}(h))<\Phi_0,\overset{\circ}{e}_0\cdot\overset{\circ}{e}_k\cdot\Phi_0
>\overset{\circ}{\ast}\overset{\circ}{e}^1.\end{eqnarray*}

Using the Clifford representation (\ref{repre}), the boundary term
is equal to (up to a constant)
\begin{equation}\beta_0(|\lambda_1|^2+|\lambda_2|^2+|\lambda_3|^2+|\lambda_4|^2)+\beta_1(-(\overline{\lambda}_2\lambda_1+\overline{\lambda}_1\lambda_2)+(\overline{\lambda}_3\lambda_4+\overline{\lambda}_4\lambda_3))$$$$+\beta_2(-\sqrt{-1}(\overline{\lambda}_2\lambda_1-\overline{\lambda}_1\lambda_2)+\sqrt{-1}(\overline{\lambda}_3\lambda_4-\overline{\lambda}_4\lambda_3))+\beta_3(|\lambda_1|^2-|\lambda_2|^2+|\lambda_4|^2-|\lambda_3|^2).\end{equation}
It can be rewritten as a quadratic form $(\bar\lambda_1,
\bar\lambda_2, \bar\lambda_3, \bar\lambda_4)Q(\lambda_1,
\lambda_2, \lambda_3, \lambda_4)^t$ where $Q$ is
\begin{equation}
\begin{pmatrix}
\beta_0+\beta_3 & -\beta_1+\sqrt{-1}\beta_2& \ & \ \\
-\beta_1-\sqrt{-1}\beta_2& \beta_0-\beta_3& \ & \ \\
\ & \ & \beta_0-\beta_3& \beta_1+\sqrt{-1}\beta_2\\
\ & \ & \beta_1-\sqrt{-1}\beta_2& \beta_0+\beta_3\end{pmatrix}.
\end{equation}

This completes the proof of the nonnegativity.

If the equality holds, then there exists at least one
non-vanishing spinor such that $\overset{\wedge}{\nabla}\phi=0$.
If $E_{\{0\}}+P_{\{0\}1} =0$, then there is $\{\phi_\alpha\}$
which forms a basis of the spinor bundle everywhere on $M$ such
that $\overset{\wedge}{\nabla}\phi_\alpha=0$. So in a local frame
$\{e_\alpha\}$ we have
$$\overline{\nabla}_i\phi_\alpha=\frac{1}{2}h_{ik}e_0\cdot e_k\cdot \phi_\alpha-\frac{\kappa}{2}e_0\cdot e_i \cdot \phi_\alpha.$$
Then
\begin{eqnarray*}\overline{\nabla}_j\overline{\nabla}_i\phi_\alpha&=&\frac{1}{2}(\overline{\nabla}_j
h_{ik})e_0\cdot e_k \cdot \phi_\alpha+\frac{1}{2}h_{ik}e_0\cdot
e_k\cdot(\frac{1}{2}h_{jl}e_0\cdot e_l-\frac{\kappa}{2}e_0\cdot
e_j)\cdot\phi_\alpha\\
&\ &-\frac{\kappa}{2}e_0\cdot e_i\cdot (\frac{1}{2}h_{jl}e_0\cdot
e_l-\frac{\kappa}{2}e_0\cdot e_j)\cdot
\phi_\alpha \\
&=&\frac{1}{2}(\overline{\nabla}_jh_{ik})e_0\cdot e_k\cdot
\phi_\alpha-\frac{1}{4}\widetilde{h}_{ik}\widetilde{h}_{jl}e_k\cdot
e_l\cdot \phi_\alpha.\end{eqnarray*} It is therefore,
\begin{eqnarray*}-\frac{1}{4}R_{ijkl}e_k\cdot e_l\cdot
\phi_\alpha&=&(\overline{\nabla}_i\overline{\nabla}_j-\overline{\nabla}_j\overline{\nabla}_i)\phi_\alpha-\overline{\nabla}_{[e_i,e_j]}\phi_\alpha\\
&=&-\frac{1}{2}(\overline{\nabla}_i\widetilde{h}_{jk}-\overline{\nabla}_j\widetilde{h}_{ik})e_0\cdot
e_k\cdot
\phi_\alpha+\frac{1}{4}(\widetilde{h}_{ik}\widetilde{h}_{jl}-\widetilde{h}_{il}\widetilde{h}_{jk})e_k\cdot
e_l\cdot \phi_\alpha\end{eqnarray*} for a basis $\{\phi_\alpha\}$.
This implies
$$\sum_{k<l}(R_{ijkl}+\widetilde{h}_{ik}\widetilde{h}_{jl}-\widetilde{h}_{il}\widetilde{h}_{jk})e_k\cdot
e_l\cdot=\sum_{k}(\overline{\nabla}_i\widetilde{h}_{jk}-\overline{\nabla}_j\widetilde{h}_{ik})e_0\cdot
e_k\cdot $$as an endomorphism of $\mathbb{S}$. Set
$$\widetilde{R}_{ijkl}=R_{ijkl}+\widetilde{h}_{ik}\widetilde{h}_{jl}-\widetilde{h}_{il}\widetilde{h}_{jk}$$
and
$$\widetilde{h}_{ijk}=\overline{\nabla}_i\widetilde{h}_{jk}-\overline{\nabla}_j\widetilde{h}_{ik}.$$
In terms of Clifford representation (\ref{repre}), we obtain
$$\begin{pmatrix} \sqrt{-1}\widetilde{R}_{ij23}&
\widetilde{R}_{ij12}+\sqrt{-1}\widetilde{R}_{ij13}&0&0\\
-\widetilde{R}_{ij12}+\sqrt{-1}\widetilde{R}_{ij13}&
-\sqrt{-1}\widetilde{R}_{ij23}&0&0\\
0&0&\sqrt{-1}\widetilde{R}_{ij23}&\widetilde{R}_{ij12}+\sqrt{-1}\widetilde{R}_{ij13}\\
0&0&-\widetilde{R}_{ij12}+\sqrt{-1}\widetilde{R}_{ij13}&-\sqrt{-1}\widetilde{R}_{ij23}\end{pmatrix}$$
$$=
\begin{pmatrix}\widetilde{h}_{ij1}&-\widetilde{h}_{ij2}-\sqrt{-1}\widetilde{h}_{ij3}&0&0\\
-\widetilde{h}_{ij2}+\sqrt{-1}\widetilde{h}_{ij3}&-\widetilde{h}_{ij1}&0&0\\
0&0&-\widetilde{h}_{ij1}&\widetilde{h}_{ij2}+\sqrt{-1}\widetilde{h}_{ij3}\\
0&0&\widetilde{h}_{ij2}-\sqrt{-1}\widetilde{h}_{ij3}&\widetilde{h}_{ij1}\end{pmatrix}.$$This
gives
$$R_{ijkl}+\widetilde{h}_{ik}\widetilde{h}_{jl}-\widetilde{h}_{il}\widetilde{h}_{jk}=0,$$
$$\overline{\nabla}_i\widetilde{h}_{jk}-\overline{\nabla}_j\widetilde{h}_{ik}=0,$$
and the theorem is proved.\qed

\mysection{The Imaginary Killing Spinors and Positive Mass Theorem
II} As mentioned in the introductory section, Maerten obtained the
positivity of a sequilinear form under the relative energy
condition. Classical linear algebra tells us that each principal
minor of this form must be nonnegative which give rise to a set of
energy-momentum inequalities (cf. Appendix \cite{Ma}). Among them,
the most interesting one might be the second order principal minor
which gives the positivity of the Lorentzian length of the mass
vector, i.e. $m_0^2-|m|^2 \geq 0$. This special inequality is
recovered in our formulism here by using the imaginary Killing
spinor.

Define the modified imaginary Killing connection as
$$\overset{\wedge}{\nabla}_i=\overline{\nabla}_i-\frac{1}{2}h_{ij}e_0\cdot e_j\cdot+\frac{\kappa\sqrt{-1}}{2}e_i\cdot.$$
Here $\bar{\nabla}$ is the Levi-Civita connection with respect to
the induced Riemannian metric on the spacelike hypersurface. The
associated hypersurface Dirac-Witten operator is
$$\overset{\wedge}{D}=\sum^3_{k=1}e_k\cdot \overset{\wedge}{\nabla}_k=D-\frac{3\kappa\sqrt{-1}}{2}.$$
The corresponding Weitzenb\"{o}ck formula is then
 \begin{eqnarray*}
\int_{\partial M}<\phi, \overset{\wedge}{\nabla}_i\phi+e_i\cdot
\overset{\wedge}{D}\phi>\ast e^i =\int_M
(|\overset{\wedge}{\nabla}\phi|^2-|\overset{\wedge}{D}\phi|^2+<\phi,\overset{\wedge}{R}\phi>)\ast
1.
 \end{eqnarray*}
Here
$$\overset{\wedge}{R}=\frac{1}{4}(Scal^{\widetilde{g}}+2\widetilde{R}_{00}+2\widetilde{R}_{0i}e_0\cdot e_i\cdot +6\kappa^2).$$(See also \cite{Ma} for
$\kappa=1$.) By the Einstein equation (\ref{ein}),
$\overset{\wedge}{R}=\frac{1}{2}(T_{00}e_0\cdot-T_{0i}e_i\cdot)e_0\cdot$
whose positivity is ensured by the standard dominant energy
condition (\ref{sdec}).

Take $\Phi_0$ as an imaginary Killing spinor along the hyperbolic
space (\ref{imk}) and extend it smoothly the whole manifold.
Consider the elliptic boundary problem $\overset{\wedge}{D}\phi=0$
with $\phi$ asymptotic to $\Phi_0$ at infinity. Using the Clifford
representation (\ref{repre}), the boundary term thus can be
written as a quadratic form $(\bar\lambda_1, \bar\lambda_2,
\bar\lambda_3, \bar\lambda_4)Q(\lambda_1, \lambda_2, \lambda_3,
\lambda_4)^t$. Here $Q$ is a Hermitian $4\times 4$ matrix
\begin{equation}\label{Q}
{\small\begin{pmatrix}
E_{\{0\}}+E_{\{3\}} & E_{\{1\}} & -P_{\{0\}2}+P_{\{3\}2} & P_{\{1\}2}- P_{\{2\}3}\\
\ & -\sqrt{-1}E_{\{2\}} & -\sqrt{-1}(P_{\{0\}3}+P_{\{3\}3}) & -\sqrt{-1}(P_{\{2\}2}-P_{\{1\}3})\\
\ &\ & \ &\ \\
E_{\{1\}} & E_{\{0\}}-E_{\{3\}} & -P_{\{1\}2}+P_{\{2\}3} & P_{\{0\}2}+P_{\{3\}2}\\
+\sqrt{-1}E_{\{2\}} & \ &-\sqrt{-1}(P_{\{2\}2}+P_{\{1\}3})& +\sqrt{-1}(P_{\{0\}3}+P_{\{3\}3})\\
\ & \ & \ & \ \\
-P_{\{0\}2}+P_{\{3\}2}& -P_{\{1\}2}+P_{\{2\}3}& E_{\{0\}}+E_{\{3\}}& -E_{\{1\}}\\
+\sqrt{-1}(P_{\{0\}3}+P_{\{3\}3})& +\sqrt{-1}(P_{\{1\}3}+P_{\{2\}2})&\ & +\sqrt{-1}E_{\{2\}} \\
\ & \ & \ & \ \\
P_{\{1\}2}-P_{\{2\}3}& P_{\{0\}2}+P_{\{3\}2}& -E_{\{1\}}& E_{\{0\}}-E_{\{3\}}\\
+\sqrt{-1}(P_{\{2\}2}-P_{\{1\}3})&
-\sqrt{-1}(P_{\{0\}3}+P_{\{3\}3})& -\sqrt{-1}E_{\{2\}} & \

\end{pmatrix}}\end{equation}

Therefore, we have reached our second positive mass theorem.
\begin{theo}\label{PET2} Let $(M,g_{ij},h_{ij})$ be a 3-dimensional asymptotically AdS initial data set which satisfies the standard dominant energy condition (\ref{sdec}).
Then $Q$ is nonnegative. Moreover, if $Q=0$, then we have the
following equations on $M$:
 \begin{eqnarray}
R_{ijkl}=(-\kappa^2)(\delta_{ik}\delta_{jl}-\delta_{il}\delta_{jk})+h_{il}h_{jk}-h_{ik}h_{jl},\label{identity1}
 \end{eqnarray}
 \begin{eqnarray}
\overline{\nabla}_ih_{jk}-\overline{\nabla}_jh_{ik}=0
\label{identity2}
 \end{eqnarray}
where $R_{ijkl}$ is the Riemann curvature tensor of $(M,g)$. These
are the Gauss and Codazzi equations of the isometric embedding in
the AdS spacetime.
\end{theo}
\rms \\
(1) The dominant energy condition in physics also implies that
$T^{00} \geq |T^{\alpha\beta}|$. If $Q=0$, then
$T^{\alpha\beta}=0$. This together with (\ref{identity1}),
(\ref{identity2}) imply that $N$ is AdS along $M$, i.e.  $N$ has
constant sectional
curvature $-\kappa^2$ along $M$.\\
(2) The energy-momentum matrix in (\ref{Q}) is different from the
one obtained in \cite{Ma} for (1+3)-dimensional spacetimes. It
should be related to the representation of spin group. However,
the method presented in our paper is consistent.

The positivity of the $2\times 2$ principal minor in (\ref{Q})
also implies the positivity of hyperbolic mass:
\begin{coro}
In particular, we have \begin{equation}E_{\{0\}} \geq
\sqrt{\sum^3_{i=1}E_{\{i\}}^2}.\end{equation}\end{coro}

Clearly, the rigidity conclusion follows from the fact that when
$Q=0$, there exists $\{\phi_\alpha\}$ which forms a basis of the
spinor bundle everywhere on $M$ such that
$\overset{\wedge}{\nabla}\phi_\alpha=0$. Maerten \cite{Ma} also
obtained this via the construction of the Killing initial data in
\cite{BC}. In addition, he discusses the isometric embedding in a
stationary pp-wave spacetime when the energy-momentum matrix is
degenerate.
\\[0.4cm]
{\bf\Large\noindent Acknowledgements}\\[0.4cm]
N. Xie is partially supported by Doctoral Foundation of Ministry
of Education of China under grant 20030246001, Eurasia-Pacific
Uninet Technologiestipendien China \& Mongolei 2005/2006
(Doktorat) and Fudan Postgraduate Students Innovation Project. He
would also thank Profs. C.H. Gu and H.S. Hu for their consistent
encouragements. X. Zhang is partially supported by National
Natural Science Foundation of China under grant 10421001,
NKBRPC(2006CB805905) and Innovation Project of Chinese Academy of
Sciences.

N. Xie\\
Institute of
Mathematics\\ School of Mathematical Sciences\\ Fudan University\\
Shanghai 200433, PR China\\ nqxie@fudan.edu.cn\\[0.5cm] X. Zhang\\ Institute of Mathematics\\ Academy of Mathematics and
System Sciences\\ Chinese Academy of Sciences\\ Beijing 100080, PR
China\\ xzhang@amss.ac.cn
\end{document}